\documentclass[11pt]{article}
\usepackage{graphicx}
\usepackage{amsmath,amsthm,latexsym,amsbsy}
\newtheorem{theorem}{Theorem}
\newtheorem{remark}{Remark}

\newtheorem{lemma}{Lemma}
\newcommand{\h}{\hspace*{.24in}}
\def\geqslant {\ge}
\def\leqslant {\le}
\def\bq{\begin{equation}}
\def\eq{\end{equation}}
\title{Determination of the body force of a two$-$dimensional isotropic elastic body}
\author{\normalsize DANG DUC TRONG$^a$, ~PHAM NGOC DINH ALAIN$^b$,\\ \normalsize PHAN THANH NAM$^a$  and TRUONG TRUNG
TUYEN$^c$\\
\\\small\it $^a$Mathematics Department, HoChiMinh City National University, Viet Nam
\\\small\it $^b$Mathematics Department, Mapmo UMR 6628, BP 67-59, 45067 Orleans cedex, France
\\\small\it $^c$Department of mathematics, Indiana University, Rawles Hall
, Bloomington, IN 47405}
\date{{}}
\begin{document}
\maketitle
\begin{abstract} Let $\Omega$ represent a two$-$dimensional isotropic elastic body. We consider the problem of determining the body force
 $F$ whose form $\varphi(t)(f_1(x),f_2(x))$ with $\varphi$ be given inexactly. The problem is nonlinear and ill-posed.
  Using the Fourier transform, the methods of Tikhonov's regularization and truncated integration,
  we construct a regularized solution from the data given inexactly and derive the explicitly error estimate.
\\\h MSC 2000: 35L20, 35R30, 42B10, 70F07, 74B05.
\\\h Key words: body force, elastic body, Fourier transform, ill$-$posed problem, Tikhonov's regularization, truncated integration.
\end{abstract}
\text{}\\
{\bf 1. Introduction}
\\\\\h Let $\Omega=(0,1)\times(0,1)$ represent a two$-$dimensional isotropic elastic body. For each $x:=(x_1,x_2)\in \Omega$, we denote by $u=(u_1(x,t)$, $u_2(x,t))$ the displacement, where $u_j$ is the displacement in the $x_j-$ direction, for all $j\in\{1,2\}$. As known, $u$ satisfies the Lam\'e system (see, e.g., \cite{TG,MHS})
$$
\frac{{\partial^2 u}}
{{\partial t^2 }} = \mu \Delta u + \left( {\lambda  +\mu } \right)\nabla \left( {div(u)} \right) + F
$$
where $F:=(F_1, F_2)$ is the body force, $div(u)=\nabla \cdot u =\partial u_1/\partial x_1+\partial u_2/\partial x_2$, and $\lambda$, $\mu$ are Lam\'e constants. We shall assume that the boundary of the elastic body is clamped and the initial conditions are given.
\\\h In this paper, we shall consider the problem of determining the body force $F$. The problem is a kind of inverse source problems. The inverse source problems are investigated in many aspects such as the uniqueness, the stability and the regularization. There are many papers devoted to the uniqueness and the stability problem. In \cite{I}, Isakov disscused the problem of finding a pair of functions $(u,f)$ satisfying
$$ cu_{tt}-\Delta u=f $$
where $f$ is independent of $t$. He proved that using some preassumptions on $f$, from the {\it final overdetermination}
$$u(x,T)=h(x)$$
, we get the uniqueness of $(u,f)$.
\\\h  As shown in \cite{MMM}, the body force (in the form $\phi(t)f(x)$) will be defined
uniquely from an observation of
surface stress (the {\it lateral overdetermination}) given on a suitable boundary of $\Omega\times
(0,T)$. In the paper, the authors also gave an abstract formula of reconstruction.
\\\h Another inverse source problem is one of finding the heat source
$F(x,t,u)$ satisfying
$$ u_t-\Delta u=F.$$
The problem was considered intensively in the last century. The problem with the final overdetermination was studied by Tikhonov in 1935 (see \cite{T}). He proved the uniqueness of problem with prescribed lateral and final data. In the last three decades, the problem is considered by many authors (see \cite{TLD,TQD,MI1,MI2,DUK,GL}). Although we have many works on the uniqueness and the stability of inverse source problems, the literature on the regularization problem is quite scarce. Very recently, in \cite{TLD,TQD} , the authors considered the regularization problem under both the lateral and the final overdetermination. The  ideas of using the Fourier transform and truncated integration in the two papers are used in the present paper.  We also consider the regularization problem under the final data and prescribed surface stress.
\\\h To get the lateral overdetermination, some mechanical arguments are in order. Let $\sigma_1,\sigma_2, \tau$ be the stresses (see \cite{TG,MHS}) defined by
$$
\begin{gathered}
  \tau  = \mu \left( {\frac{{\partial u_1 }}
{{\partial x_2 }} + \frac{{\partial u_2 }}
{{\partial x_1 }}} \right) \hfill \\
  \sigma _j  = \lambda div(u)+ 2\mu \frac{{\partial u_j }}
{{\partial x_j }},\h j\in\{1,2\} \hfill \\
\end{gathered}
$$
We shall assume that the surface stress is given on the boundary of the body, i.e.,
\[
\left( \begin{gathered}
  \sigma _1 \h \tau \hfill \\
  \tau\h\sigma _{\text{2}}  \hfill \\
 \end{gathered}  \right)\left( \begin{gathered}
  n_1  \hfill \\
  n_2  \hfill \\
 \end{gathered}  \right) = \left( \begin{gathered}
  X_1  \hfill \\
  X_2  \hfill \\
 \end{gathered}  \right)
\]
where  $X=(X_1,X_2)$ is given on $\partial\Omega$, and $n=(n_1,n_2)$ is the outward unit normal vector of $\partial \Omega$.
\\\h As discussed, our problem is severely ill-posed. Hence, to simplify the problem, a preassumption on the form
of the body force is needed. We shall use the separable form force as in \cite{MMM}
$$
(F_1(x,t),F_2(x,t))= \varphi(t)(f_1(x),f_2(x))
$$
where $\varphi$ is given inexactly. The form is issued from an approximated model for elastic wave generated from a point dislocation source (see, e.g., \cite{MMM,AKR}). But, since $\varphi$ is inexact, our problem is nonlinear. Morever, the problem is still ill-posed because the measured data is not only inexact but also non-smooth.
\\\h Precisely, we consider the problem of identifying a pair of functions  $(u,f)$ satisfying the system:
\bq
  \frac{{\partial ^2 u_j }}
{{\partial t^2 }} = \mu \Delta u_j  + (\lambda  + \mu )\frac{\partial  }
{{\partial x_j }}div(u)+ \varphi (t)f_j (x) ,\forall j\in\{1,2\} \label{1}
\eq
for $(x,t)\in \Omega\times (0,T)$, where $\mu,\lambda$ are real constants satisfying $\mu>0$ and $\lambda+2\mu>0$.
\\Since the boundary of the elastic body is clamped, the displacement $u=(u_1,u_2)$ satisfies the boundary condition
\bq (u_1(x,t),u_2(x,t))=(0,0),\h x\in \partial \Omega \label{bcd}\eq
In addition, the initial and final displacement are given in $\Omega$
\bq
\left\{ \begin{gathered}
   (u_1(x,0),u_2(x,0))=(u_{01}(x),u_{02}(x))\hfill \\
   \left( {\frac{\partial u_1}{\partial t}(x,0),\frac{\partial u_2}{\partial t}(x,0)}\right)
=(u_{01}^*(x),u_{02}^*(x))\hfill \\
   (u_1(x,T),u_2(x,T))=(u_{T1}(x),u_{T2}(x))\hfill \\
\end{gathered}  \right.
 \eq
Finally, the surface stress is given on $\partial \Omega$
\
\bq
 \left\{ \begin{gathered}
  n_1 \sigma _1  + n_2 \tau = {X_1 }  \hfill \\
  n_2 \sigma _2  + n_1 \tau = {X_2 }  \hfill \\
\end{gathered} \right.
\eq
\h We shall assume that the data of the system $(1)-(4)$
\[
I = (\varphi ,X,u_0 ,u_0^*,u_T )\in
\left( {L^1(0,T),(L^1(0,T,L^1(\partial \Omega)))^2, (L^1(\Omega))^2,(L^1(\Omega))^2,(L^1(\Omega))^2} \right)
\]
are given inexactly since they are results of experimental measurements. The system $(1)-(4)$ usually has no solution; moreover, even if the solution exists, it does not  depend continously on the given data. Hence, a regularization is in order. Denoting by $I_{ex}$ the exact data, which are probably unknown, corresponding to an exact solution $(u_{ex},f_{ex})$ of the system $(1)-(4)$ , from the inexact data $I_{\varepsilon}$ approximating $I_{ex}$, we shall construct a regularized solution $f_{\varepsilon}$ approximating $f_{ex}$ .
\\\h In fact, using the Fourier transform, we shall reduce our problem to finding the solutions of the binomial equations whose binomial term is an entire function (see Lemma 1). In this case, the problem is unstable in the neighborhood of zeros of the entire function. The zeroes can be seen as singular values. Using the  method of Tikhonov's regularization and truncated integration, we shall eliminate the singular values to regularize our problem. Error estimates are given.
\\\h The remainder of the paper is divided into two sections. In Section 2, we shall set some  notations and state our main results. In Section 3, we give the proofs of the results.
\\\\{\bf 2. Notations and main results}
\\\\\h We recall that $\Omega=(0,1)\times (0,1)$. We always assume that the data $I = (\varphi ,X,u_0 ,u_T ,u_T^* )$  belong to
$$
\left( {L^1(0,T),(L^1(0,T,L^1(\partial \Omega)))^2, (L^1(\Omega))^2,(L^1(\Omega))^2,(L^1(\Omega))^2} \right)
$$
For all $\xi  = (\xi _1 ,\xi _2 ),\zeta  = (\zeta _1 ,\zeta _2 )\in R^2$, we set $\xi \cdot\zeta=\xi _1 \zeta _1  + \xi _2 \zeta _2  $ and $\left| \xi  \right| = \sqrt {\xi \cdot\xi }$.
\\\h We first have the following lemma.
\begin{lemma} If $u\in (C^2([0,T];L^2(\Omega))\cap L^2(0,T;H^2(\Omega)))^2$, $f\in (L^2(\Omega))^2$ satisfy $(1)-(4)$ corresponding the data I, then for all $\alpha=(\alpha_1,\alpha_2)\in R^2 \backslash \{ 0\} $, we have
\[
2D(I)  .\int\limits_\Omega  {f_j (x).\cos (\alpha  \cdot x)dx}  =g_j(I),\h\forall j\in\{1,2\}
\]
where
$$D(I)=D_1(I).D_2(I), g_j(I)=\frac{2}{|\alpha|^2}(\alpha _j D_2(I) h_0  + D_1(I) h_j)$$ with
$$
D_1(I)  = \int\limits_0^T {\varphi (T-t)\sin (\sqrt {\lambda  + 2\mu } \left| \alpha  \right|t)dt}, D_2(I)  = \int\limits_0^T  {\varphi (T-t)\sin (\sqrt \mu  \left| \alpha  \right|t)dt}
$$
$$
\begin{gathered}
 h_0(I)=-\sin (\sqrt {\lambda  + 2\mu } |\alpha |T).\int\limits_\Omega  {(\alpha  \cdot u_0^* ).\cos (\alpha  \cdot x)dx}\hfill\\
\h +\sqrt {\lambda  + 2\mu }.|\alpha|. \int\limits_\Omega  {(\alpha  \cdot u_T ).\cos (\alpha  \cdot x)dx}  \hfill \\
  \h- \sqrt {\lambda  + 2\mu }.|\alpha|.\cos (\sqrt {\lambda  + 2\mu } |\alpha |T).\int\limits_\Omega  {(\alpha  \cdot u_0 ).\cos (\alpha  \cdot x)dx} \hfill\\
  \h - \int\limits_0^T {\int\limits_{\partial \Omega } {\sin (\sqrt {\lambda  + 2\mu } |\alpha |(T-t))(\alpha  \cdot X).\cos (\alpha  \cdot x) d\omega dt} }  \hfill \\
\end{gathered}
$$
$$
\begin{gathered}
h_j(I)
   = -\sin (\sqrt \mu  |\alpha |T).\int\limits_\Omega  {(|\alpha |^2u_{0j}^*  - \alpha _j (\alpha  \cdot u_0^* )).\cos (\alpha  \cdot x)dx}  \hfill\\
\h+ \sqrt \mu .|\alpha| .\int\limits_\Omega  {(|\alpha |^2   u_{Tj}  - \alpha _j (\alpha  \cdot u_T )).\cos (\alpha  \cdot x)dx}  \hfill \\
   \h- \sqrt \mu.|\alpha|.  \cos (\sqrt \mu  |\alpha |T).\int\limits_\Omega  {(|\alpha |^2  u_{0j}  - \alpha _j (\alpha  \cdot u_0)).\cos (\alpha  \cdot x)dx}  \hfill \\
 \h  - \int\limits_0^T {\int\limits_{\partial \Omega } {\sin (\sqrt \mu  |\alpha |(T-t))(|\alpha |^2 X_j  - \alpha _j (\alpha  \cdot X)).\cos (\alpha  \cdot x) d\omega dt} },\forall j\in\{1,2\}. \hfill \\
\end{gathered}
$$
\end{lemma}
From Lemma 1, we consider the function
\[
D(I)  = \int\limits_0^T {\varphi (T-t)\sin (\sqrt {\lambda  + 2\mu } \left| \alpha  \right|t)dt} .\int\limits_0^T {\varphi (T-t)\sin (\sqrt \mu  \left| \alpha  \right|t)dt}
\]
The problem is unstable in the neighborhood of zeros of this function. However, from the properties of analytic function, we can show that if $\varphi\not \equiv 0$ then this function differ from $0$ for almost every where in $R^3$. Furthermore, using the idea of Theorem 4 in \cite{TT}, we get the following lemma.
\begin{lemma} Let $\tau,q$ be positive constants, $\varphi_0\in L^1(0,T)\backslash \{0\}$ and $D(\varphi_0,\tau):R^2\to R$
\[
D(\varphi_0,\tau)(\alpha ) = \int\limits_0^T {\varphi_0 (t)\sin (\sqrt {\tau}|\alpha| t)dt}
\]
Then $D(\varphi_0,\tau)\ne 0$ for a.e $\alpha\in R^2$. Moreover, if we put
\[
R_\varepsilon=\frac{q}
{{9eT}}.\frac{{\ln (\varepsilon ^{ - 1} )}}
{{\ln (\ln (\varepsilon ^{ - 1} ))}},\h \forall \varepsilon >0
\]
then the Lebesgue measure of the set
\[
B(\varphi_0,\tau,\varepsilon)=\left\{ {\alpha\in B(0,R_\varepsilon)  ,|D(\varphi_0,\tau)(\alpha )| \leqslant \varepsilon ^q} \right\}
\]
is less than $R_\varepsilon^{-1}$ for $\varepsilon>0$ small enough, where $B(0,R_\varepsilon)$ is the open ball in $R^2$ .
\end{lemma}
Lemma 1 and Lemma 2 imply immediately the uniqueness result.
\begin{theorem} Let $u,u^*\in (C^2([0,T];L^2(\Omega)) \cap L^2(0,T;H^2(\Omega)))^2$, $f,f^*\in (L^2(\Omega))^2$. If $(u,f)$, $(u^*,f^*)$ satisfy $(1)-(4)$ corresponding the same data $I$, and $\varphi\not \equiv 0$, then
$$(u,f)=(u^*,f^*)$$
\end{theorem}
Let $(u_{ex},f_{ex})$ be the exact solution of the system $(1)-(4)$ corresponding the exact data $I_{ex}=(\varphi_{ex} ,X_{ex} ,u_0^{ex} ,u_0^{*ex} ,u_T^{ex} )$. Notice that, if we assume
\bq
\begin{gathered}
   u_{ex}\in (C^2([0,T];L^2(\Omega))\cap L^2([0,T];H^2(\Omega)))^2, f_{ex}\in (L^2(\Omega))^2,\varphi_{ex}\in L^1(0,T)\backslash \{0\} \\
\end{gathered} \label{H}
\eq
then for all $j\in\{1,2\}$,
\[
F(\widetilde f_{jex})(\alpha ) =2\int\limits_\Omega  {f_{jex} (x)\cos (\alpha  \cdot x)dx}
=
 \frac{{g_{j} (I_{ex} )}}
{{D(I_{ex} )}}
\]
for  a.e $\alpha\in R^2$, where $g_j$, $D$ are defined by Lemma 1, $\widetilde f_{jex}:R^2\to R$ is defined by $\widetilde f_{jex}(x)=\chi(\Omega)f_{jex}(x)+\chi(-\Omega)f_{jex}(-x)$, and  $F$ is the Fourier transform in $R^2$.
\\\h From approximate data $I_\varepsilon = (\varphi ,X,u_0 ,u_0^*,u_T  )$ satisfying
\bq
\begin{gathered}
 \left\| {\varphi  - \varphi_{ex} } \right\|_{L^1 (0,T)}  \leqslant \varepsilon ,\left\| {X_j  - X_j^{ex} } \right\|_{L^1 (0,T,L^1(\partial \Omega))}  \leqslant \varepsilon ,\left\| {u_{0j}  - u_{0j}^{ex} } \right\|_{L^1 (\Omega )}  \leqslant \varepsilon  \hfill \\
  \left\| {u_{0j}^*  - u_{0i}^{*ex} } \right\|_{L^1(\Omega )}  \leqslant \varepsilon ,\left\| {u_{Tj}  - u_{Tj}^{ex} } \right\|_{L^1(\Omega )}  \leqslant \varepsilon,\h\forall j\in \{1,2\}  \hfill \\
\end{gathered} \label{xx}
\eq
, we construct a regularized solution $f_{\varepsilon}=(f_{1\varepsilon},f_{2\varepsilon})$ whose Fourier transform is
\[
F( f_{j\varepsilon})(\alpha ) =\chi(B(0,R_{\varepsilon})).\frac{{g_{j} (I_{\varepsilon} ).D(I_\varepsilon)}}
{{\delta _\varepsilon+\left( {D(I_\varepsilon  )} \right)^2}},\forall \alpha\in R^2\backslash \{0\}
\]
where
\bq
q=\frac{1}{7}, \delta _\varepsilon=\varepsilon^{\frac{1+6q}{2}}, R_{\varepsilon}=\frac{q}{9eT}.\frac{ln(\varepsilon^{-1})}{ln(ln(\varepsilon^{-1}))}\label{ths}
\eq
We have two regularization results.
\begin{theorem} Let $(u_{ex},f_{ex})$ be the exact solution of the system $(1)-(4)$ corresponding the exact data $I_{ex}$, and $(\ref{H})$ hold. Then from the given data $I_\varepsilon$  satisfying $(\ref{xx})$, we can construct a regularized solution $f_{\varepsilon}\in (C(\overline\Omega))^2$ such that
\[
\mathop {\lim }\limits_{\varepsilon  \to 0} \left\| {f_{j\varepsilon }  - f_{jex} } \right\|_{L^2 (\Omega )}  = 0,\h\forall j\in \{1,2\}
\]
If we assume, in addition, that $f_{ex}\in (H^1(\Omega))^2$,then
\[
\left\| {f_{j\varepsilon }  - f_{jex} } \right\|^2_{L^2 (\Omega )}  \leqslant 63eT\left( {66\left\| {f_{jex} } \right\|_{H^1 (\Omega )}^2  + (2\pi) ^{ -2} } \right).\frac{{\ln (\ln (\varepsilon ^{ - 1} ))}}
{{\ln (\varepsilon ^{ - 1} )}},\h \forall j\in \{1,2\}
\]
for $\varepsilon>0$ small enough.
\end{theorem}
\begin{theorem} Let $(u_{ex},f_{ex})$ be the exact solution of the system $(1)-(4)$ corresponding the exact data $I_{ex}$, and $(\ref{H})$ hold. We assume, in addition, that
\[
\int\limits_{R^2 } {\left| {\int\limits_\Omega  {f_{jex} (x).\cos (\alpha  \cdot x)dx} } \right|d\alpha  < \infty },\h\forall j\in\{1,2\}
\]
Then from the given data $I_\varepsilon$  satisfying $(\ref{xx})$, we can construct a regularized solution $f_{\varepsilon}\in (C(\overline\Omega))^2$, which coincides the one in Theorem 2, such that
$$
\mathop {\lim }\limits_{\varepsilon  \to 0} \left\| {f_{j\varepsilon }  - f_{jex} } \right\|_{L^\infty  (\Omega )}  = 0,\h\forall j\in \{1,2\}
$$
\end{theorem}
\text{}
\\{\bf 3. Proofs of the results}
\\\\{\bf Proof of Lemma 1}
\begin{proof} Let $\alpha=(\alpha_1,\alpha_2)\in R^2$ and $G=\cos(\alpha\cdot x)$. Notice that  the j$-$th equation of the system (\ref{1}) can rewrite
\[
\frac{{\partial^2 u_j }}
{{\partial t^2 }} = \frac{{\partial \sigma _j }}
{{\partial x_j }} + \frac{{\partial \tau }}
{{\partial x_k }} + \varphi (t)f_j (x),\h \{j,k\}=\{1,2\}
\]
Getting the inner product (in $L^2(\Omega)$) of the equation and $G$ and using the condition $(\ref{bcd})$, for $\{j,k\}=\{1,2\}$, we get
\bq
\begin{gathered}
  \frac{d}
{{dt^2 }}\int\limits_\Omega  {u_j G}  = \int\limits_{\partial \Omega } {(n_j \sigma _j  + n_k \tau )Gd\omega }  - \int\limits_\Omega  {\sigma _j \frac{{\partial G}}
{{\partial x_j }}dx}  - \int\limits_\Omega  {\tau \frac{{\partial G}}
{{\partial x_k }}dx}  + \varphi (t)\int\limits_\Omega  {f_j Gdx}  \hfill \\
   = \int\limits_{\partial \Omega } {X_j Gd\omega }  - \mu \left| \alpha  \right|^2 \int\limits_\Omega  {u_j Gdx}  - (\lambda  + \mu )\alpha_j\int\limits_\Omega  {(\alpha  \cdot u)Gdx}  + \varphi (t)\int\limits_\Omega  {f_j Gdx}  \hfill \\
\end{gathered}  \label{pt0}
\eq
Multiplying $(\ref{pt0})$ by $\alpha_j$, then getting the sum for $j=1,2$, we obtain
\bq
  \frac{d}
{{dt^2 }}\int\limits_\Omega  {(\alpha\cdot u)G dx} = \int\limits_{\partial \Omega } {(\alpha \cdot X) Gd\omega }
   - (\lambda  + 2\mu )|\alpha|^2\int\limits_\Omega  {(\alpha\cdot u)G dx}  + \varphi (t)\int\limits_\Omega  {(\alpha\cdot f)G dx}
\label{pt1}
\eq
Multiplying $(\ref{pt0})$ by $|\alpha|^2$ and multiplying $(\ref{pt1})$ by $-\alpha_j$, then getting the sum of them, we have
\bq
\begin{gathered}
  \frac{d}
{{dt^2 }}\int\limits_\Omega  {\left( {\left| \alpha  \right|^2 u_j  - \alpha _j .\left( {\alpha  \cdot u} \right)} \right)Gdx}  = \int\limits_{\partial \Omega } {\left( {\left| \alpha  \right|^2 X_j  - \alpha _j .\left( {\alpha  \cdot X} \right)} \right)G dx}  \hfill \\
\h   - \mu \left| \alpha  \right|^2 \int\limits_\Omega  {\left( {\left| \alpha  \right|^2 u_j  - \alpha _j .\left( {\alpha  \cdot u} \right)} \right)Gdx}  + \varphi(t)\int\limits_\Omega  {\left( {\left| \alpha  \right|^2 f_j  - \alpha _j .\left( {\alpha  \cdot f} \right)} \right)Gdx}  \hfill \\
\end{gathered} \label{pt2}
\eq
\h We consider $(\ref{pt1})$ and $(\ref{pt2})$ as the differential equations whose form
\bq y''+\eta^2y=h(t)\label{ptr}\eq
where $\eta$ is a real constant and $y(0)$, $y'(0)$, $y(T)$ are given. Getting the inner product (in $L^2(0,T)$) of $(\ref{ptr})$ and $\sin(\eta (T-t))$, we have
\bq
-y'(0)sin(\eta T)+\eta y(T)-\eta y(0)cos(\eta T)=\int\limits_0^T {h(T-t)\sin (\eta t)dt} \label{nt}
\eq
Applying $(\ref{nt})$ to $(\ref{pt1})$ with $\eta=\sqrt{(\lambda+2\mu)}|\alpha|$ and $y = \int\limits_\Omega  {(\alpha \cdot u).Gdx}$, we get
\bq
D_1(I) .\int\limits_\Omega  {(\alpha  \cdot f).Gdx}=h_0(I) \label{pt3}
\eq
where $D_1(I),h_0(I)$ are defined by Lemma 1.
\\Similarly, applying $(\ref{nt})$ to $(\ref{pt2})$ with $\eta=\sqrt{\mu}|\alpha|$ and $y = \int\limits_\Omega  {(|\alpha|^2u_j-\alpha_j.(\alpha \cdot u)).Gdx}$, we get
\bq
  D_2(I).\int\limits_\Omega  {(|\alpha |^2 f_j  - \alpha _j (\alpha  \cdot f)).Gdx} =h_j(I),\h\forall j\in\{1,2\} \label{pt4}
\eq
where $D_2(I),h_j(I)$ are defined by Lemma 1.
\\\h Multiplying $(\ref{pt3})$ by $\alpha_jD_2(I)$ and multiplying $(\ref{pt4})$ by $D_1(I)$, then getting the sum of them, we obtain the result of Lemma 1.
\end{proof}
\text{}\\{\bf Proof of Lemma 2}
\begin{proof} Put $\widetilde\varphi_0:R\to R$
\[
\widetilde\varphi_0 (t) = \frac{1}
{2}\left\{ \begin{gathered}
  \varphi_0 (t)\h t\in(0,T) \hfill \\
   - \varphi_0 ( - t)\h t\in (-T,0) \hfill \\
  0\h t\notin (-T,T) \hfill \\
\end{gathered}  \right.
\]
and $\phi:C \to C$
\[
\phi(z) = \int\limits_{ - \infty }^\infty  {e^{-itz} \widetilde\varphi_0 (t)dt}  = \int\limits_{ - T}^T {e^{-itz} \widetilde\varphi_0 (t)dt}
\]
Then $\phi$ is an entire function and $D(\varphi_0,\tau)(\alpha)=i\phi( \sqrt {\tau}|\alpha|)$. Because $\widetilde\varphi_0\not\equiv 0$, its Fourier transform (in R) does not coincide 0. Therefore, there exists $z_0\in R$ such that $|\phi(z_0)|=C_1>0$. Thus $\phi\not\equiv 0$. Since $\phi$ is an entire function, its zeros set is either finite or countable. Consequently, $D(\varphi_0,\tau)(\alpha)\ne 0$ for a.e $\alpha\in R^2$.
\\\h To estimate the measure of $B(\varphi_0,\tau ,\varepsilon)$, we shall use the following result (see Theorem 4 of $\$ 11.3$ in \cite{Ya}).
\begin{lemma} Let $f(z)$ be a function analytic in the disk $\{z:|z|\le 2eR\}$, $|f(0)|=1$, and let $\eta$ be an arbitrary small positive number. Then the estimate
\[
\ln |f(z)| >  - \ln (\frac{{15e^3 }}
{\eta }).\ln (M_f (2eR))
\]
is valid everywhere in the disk $\{z:|z|\le R\}$ except a set of disks $(C_j)$ with sum of radii $\sum r_j\le \eta R$. Where $M_f (r) = \mathop {\max }\limits_{|z| = r} |f(z)|$.\end{lemma}
\h Returning Lemma 2, we put $\phi_1:C\to C$
\[
\phi _1 (z) = \frac{{\phi (z + z_0)}}
{{C_1}}
\]
Then $\phi_1$ is an entire function, $\phi_1(0)=1$, and for all $z \in C,\left| z \right| \leqslant 2eR$,
\[
C_1\left| {\phi _1 (z)} \right| =\left| {\int\limits_{ - T}^T {e^{ - it(z + z_0 )} \widetilde\varphi_0(t)} } \right| \leqslant e^{2eRT}.\int\limits_{-T}^T {\left| {\widetilde\varphi_0 (t)} \right|dt}
 = e^{2eRT}\left\| {\varphi_0} \right\|_{L^1 (0,T)}
\]
For $\varepsilon>0$ small enough, applying Lemma 3 with $R=\frac{4}{3}R_\varepsilon$ and $\eta  = \frac{{\sqrt \tau  }}{{8\pi R_\varepsilon ^3}}$, we get
\[
\begin{gathered}
  \ln \left| {\phi _1 (z)} \right| >  - \left[ {3\ln R_\varepsilon   + \ln (\frac{{8\pi }}
{{\sqrt \tau  }}) + \ln (15e^3 )} \right].\left[ {\frac{8}{3}.eTR_\varepsilon   + \ln \left( {\frac{{\left\| \varphi_0  \right\|_{L^1 (0,T)} }}
{{C_1}}} \right)} \right]\hfill\\
  \h\h\h > - \frac{17}{2}T.R_\varepsilon  \ln R_\varepsilon >-q\ln(\varepsilon^{-1})-\ln(C_1)= \ln (\frac{\varepsilon ^q}{C_1})\hfill \\
\end{gathered}
\]
for all $|z|\le \frac{4}{3}R_\varepsilon$ except a set of disks $\{ B(z_j,r_j)\}_{j \in J}$ with sum of radii $\sum {r_i }  \leqslant \eta R = \frac{{\sqrt \tau  }}
{{6\pi R_\varepsilon ^2 }}$.
\\\h Consequently, for $\varepsilon>0$ small enough, we have $|z_0|<\frac{1}{3}R_{\varepsilon}$ and $\left| {\phi (z)} \right|= C_1 .\left| {\phi _1 (z - z_0 )} \right|\geqslant \varepsilon ^q$ for all $|z|\le R_\varepsilon$ except the set $\mathop  \cup \limits_{j \in J} B(z_j+z_0 ,r_j )$ . Hence, $B(\varphi_0, \tau, \varepsilon)$ is contained in the set $\mathop  \cup \limits_{j \in J} B_j$, where
\[
B_j=\{\alpha\in B(0,R_\varepsilon)  ,\left| {\sqrt {\tau}|\alpha|  - y_j } \right| \leqslant r_j \}
\]
with $y_j=Re(z_j+z_0)$.
\\\h If $y_j>\sqrt{\tau}R_\varepsilon+r_j$ then $B_j=\emptyset$. If $y_j\le r_j$ then $B_j  \subset B (0,\frac{2r_j}{\sqrt{\tau}} )$, so $m(B_j ) \leqslant \frac{{4\pi r_j^2 }}{\tau}$. If $r_j<y_j\le \sqrt{\tau}R_\varepsilon+r_j$ then
\[
B_j  \subset B(0,\frac{{y_j  + r_j }}
{{\sqrt \tau  }})\backslash B (0,\frac{{y_j  - r_j }}
{{\sqrt \tau  }})
\]
hence
\[
m(B_j ) \leqslant \frac{{\pi (y_j  + r_j )^2 }}
{\tau } - \frac{{\pi (y_j  -r_j )^2 }}
{\tau } = \frac{{4\pi y_j r_j }}
{\tau } \leqslant \frac{{4\pi (\sqrt \tau  R_\varepsilon   + r_j )r_j }}
{\tau }
\]
Thus we get
\[
\begin{gathered}
  m(B(\varphi, \tau,\varepsilon)) \leqslant \sum {\frac{{4\pi (\sqrt \tau  R_\varepsilon   + r_j )r_j }}
{\tau }}  + \sum {\frac{{4\pi r_j^2 }}
{\tau }}  \hfill \\
   \leqslant \frac{{4\pi R_\varepsilon  }}
{{\sqrt \tau  }}\sum {r_j }  + \frac{{8\pi }}
{\tau }(\sum {r_j } )^2  \leqslant \frac{{4\pi R_\varepsilon  }}
{{\sqrt \tau  }}.\frac{{\sqrt \tau  }}
{{6\pi R_\varepsilon ^2 }} + \frac{{8\pi }}
{\tau }.(\frac{{\sqrt \tau  }}
{{6\pi R_\varepsilon ^2 }})^2<\frac{1}
{{R_\varepsilon  }} \hfill \\
\end{gathered}
\]
for $\varepsilon>0$ small enough. The proof of Lemma 2 is completed.
\end{proof}
\text{}\\{\bf Proof of theorem 1}
\begin{proof}
Put $w=u-u^*$ and $v=f-f^*$ then $(w,v)$ satisfies $(1)-(4)$ corresponding the data
$$I=(\varphi,(0,0),(0,0),(0,0),(0,0))$$
Let $\widetilde v_j:R^2\to R$ be defined by $\widetilde v_j(x)=\chi(\Omega)  v_j(x)+\chi(-\Omega)v_j( - x)$. Lemma 1 implies that, for all $j\in\{1,2\}$, for all $\alpha\in R^2 \backslash \{ 0\} $, we get
\[
D(I).F(\widetilde v_j)(\alpha ) = 2D(I).\int\limits_\Omega  {v_j(x)\cos (\alpha  \cdot x)dx}=g_j(I)=0
\]
Applying Lemma 2 with $\varphi_0(t)=\varphi(T-t)$, we get $D(I)\ne 0$ for  a.e $\alpha\in R^2$. Therefore, $F(\widetilde v_j)\equiv 0$, and it implies that $\widetilde v_j \equiv 0$. Thus $v\equiv (0,0)$. Hence, $w$ satisfies that
\bq
\frac{{\partial^2 w}}
{{\partial t^2 }} = \mu \Delta w + \left( {\lambda  +\mu } \right)\nabla \left( {div(w)} \right)\label{eq}
\eq
Getting the inner product (in $(L^2(\Omega))^2$) of $(\ref{eq})$ and $\partial w/\partial t$, we have
\[
\frac{1}
{2}.\frac{d}
{{dt}}\sum\limits_{j = 1}^2 {\left\| {\frac{{\partial w_j }}
{{\partial t}}} \right\|_{L^2 (\Omega )}^2 }  =  - \frac{\mu }
{2}.\frac{d}
{{dt}}\sum\limits_{j = 1}^2 {\left\| {\nabla w_j } \right\|_{L^2 (\Omega )}^2 }  - \frac{{\lambda  + \mu }}
{2}.\frac{d}
{{dt}}\left\| {div(w)} \right\|_{L^2 (\Omega )}^2
\]
Integrating this equality in $(0,t)$, we get
\bq
\sum\limits_{j = 1}^2 {\left\| {\frac{{\partial w_j }}
{{\partial t}}} \right\|_{L^2 (\Omega )}^2 }  + \mu \sum\limits_{j = 1}^2 {\left\| {\nabla w_j } \right\|_{L^2 (\Omega )}^2 }  + (\lambda  + \mu )\left\| {div(w)} \right\|_{L^2 (\Omega )}^2  = 0\label{dgtn}
\eq
for all $t\in (0,T)$. Using the condition $(\ref{bcd})$, we have
\[
\begin{gathered}
  \left\| {div(w)} \right\|_{L^2 (\Omega )}^2  = \sum\limits_{j = 1}^2 {\left\| {\frac{{\partial w_j }}
{{\partial x_j }}} \right\|_{L^2 (\Omega )}^2 }  + 2{\int\limits_\Omega  {\frac{{\partial w_1 }}
{{\partial x_1 }}.\frac{{\partial w_2 }}
{{\partial x_2 }}} }
   = \sum\limits_{j = 1}^2 {\left\| {\frac{{\partial w_j }}
{{\partial x_j }}} \right\|_{L^2 (\Omega )}^2 }  + 2{\int\limits_\Omega  {\frac{{\partial w_1 }}
{{\partial x_2 }}.\frac{{\partial w_2 }}
{{\partial x_1 }}} }  \hfill \\
 \h\h\h\h  \leqslant \sum\limits_{j = 1}^2 {\left\| {\frac{{\partial w_j }}
{{\partial x_j }}} \right\|_{L^2 (\Omega )}^2 }  + {\left( {\left\| {\frac{{\partial w_1 }}
{{\partial x_2 }}} \right\|_{L^2 (\Omega )}^2  + \left\| {\frac{{\partial w_2 }}
{{\partial x_1 }}} \right\|_{L^2 (\Omega )}^2 } \right)}  = \sum\limits_{j = 1}^2 {\left\| {\nabla w_j } \right\|_{L^2 (\Omega )}^2 }  \hfill \\
\end{gathered}
\]
Since $\mu>0$ and $\lambda  + 2\mu>0$, the above inequality implies that
\[
\mu \sum\limits_{j = 1}^2 {\left\| {\nabla w_j } \right\|_{L^2 (\Omega )}^2 }  + (\lambda  + \mu )\left\| {div(w)} \right\|_{L^2 (\Omega )}^2  \geqslant 0
\]
From $(\ref{dgtn})$, we obtain $\partial w/\partial t=(0,0)$. Since $w(x,0)=(0,0)$, the proof is completed.
\end{proof}
To prove two main regularization results, we state and prove some preliminary lemmas.
\begin{lemma} Let $(u_{ex},f_{ex})$ be the exact solution of $(1)-(4)$ corresponding the exact data $I_{ex}$ satisfying $(\ref{H})$, and the given data $I_\varepsilon$  satisfying $(\ref{xx})$. Using notations of $(\ref{ths})$, we put
\[
  G_j (I_\varepsilon  ) = \chi (B(0,R_\varepsilon  )).\frac{{g_j (I_\varepsilon  )D(I_\varepsilon  )}}
{{\delta _\varepsilon   + \left( {D(I_\varepsilon  )} \right)^2 }}
\]
Then for all $j\in\{1,2\}$, we have $G_j (I_\varepsilon  ) \in L^1 (R^2 ) \cap L^2 (R^2)$; moreover, there exists a constant $C_0$ depend only on $I_{ex}$ such that for all $\varepsilon\in (0,e^{-e})$,
\[
\begin{gathered}
  \left| {G_j (I_\varepsilon  ) - F(\widetilde f_{jex} )} \right| \leqslant \chi (B(0,R_\varepsilon  ))C_0 R_\varepsilon  \varepsilon ^{\frac{{1 - 6q}}
{2}}  \hfill \\
 \h  + 2\chi (B_\varepsilon  )\left\| {f_{jex} } \right\|_{L^2 (\Omega )}  + \chi (R^2 \backslash B(0,R_\varepsilon  ))\left| {F(\widetilde f_{jex} )} \right| \hfill \\
\end{gathered}
\]
where $B_\varepsilon   = \left\{ {\alpha  \in B(0,R_\varepsilon  ),\left| {D(I_{ex} )(\alpha )} \right| \leqslant \varepsilon ^{2q} } \right\}$.
\end{lemma}
\begin{proof} First, we show that there exists a constant $C_2>0$ depend only on $I_{ex}$ such that for all $\varepsilon\in (0,e^{-e})$, $r>r_0=q/(9T)$, $j\in\{1,2\}$,
\[
\begin{gathered}
  \left\| {D(I_{ex} )} \right\|_{L^\infty  (R^2 )}  \leqslant C_2 ,\left\| {D(I_\varepsilon  ) - D(I_{ex} )} \right\|_{L^\infty  (R^2 )}  \leqslant C_2 \varepsilon  \hfill \\
  \left\| {g_j (I_{ex} )} \right\|_{L^\infty  (B(0,r))}  \leqslant C_2 r,\left\| {g_j (I_\varepsilon  ) - g_j (I_{ex} )} \right\|_{L^\infty  (B(0,r))}  \leqslant C_2 r\varepsilon  \hfill \\
\end{gathered}
\]
Recall that $D_1(I), D_2(I), h_0(I),h_j(I)$ are defined by Lemma 1. For all $\alpha \in R^3$ we have
\[
\left| {D_k (I_{ex} )} \right| \leqslant \left\| {\varphi _{ex} } \right\|_{L^1 (0,T)} ,\left| {D_k (I_\varepsilon  ) - D_k (I_{ex} )} \right| \leqslant \left\| {\varphi _\varepsilon   - \varphi _{ex} } \right\|_{L^1 (0,T)}  \leqslant \varepsilon
\]
for all $k\in\{1,2\}$. Hence, $| {D(I_{ex} )} |\leqslant \left\| {\varphi _{ex} } \right\|_{L^1 (0,T)}^2$ and
\[
\begin{gathered}
  \left| {D(I_\varepsilon  ) - D(I_{ex} )} \right| = \left| {D_1 (I_\varepsilon  ) - D_1 (I_{ex} )} \right|.\left| {D_2 (I_\varepsilon  )} \right| + \left| {D_1 (I_{ex} )} \right|.\left| {D_2 (I_\varepsilon  ) - D_2 (I_{ex} )} \right| \hfill \\
   \leqslant \varepsilon .(\left\| {\varphi _{ex} } \right\|_{L^1 (0,T)}  + \varepsilon ) + \left\| {\varphi _{ex} } \right\|_{L^1 (0,T)} .\varepsilon  \leqslant (2\left\| {\varphi _{ex} } \right\|_{L^1 (0,T)}+e^{-e}).\varepsilon  \hfill \\
\end{gathered}
\]
A straightforward calculation show that, for all $\alpha\in B(0,r)\backslash \{0\}$, we have
\[
\begin{gathered}
  \left| {\alpha_j h_0 (I_{ex} )} \right| \leqslant C_3 r\left| \alpha  \right|^2 ,\left| {\alpha_j (h_0 (I_\varepsilon  ) - h_0 (I_{ex} ))} \right| \leqslant C_3 r\left| \alpha  \right|^2 \varepsilon , \hfill \\
  \left| {h_j (I_{ex} )} \right| \leqslant C_3 r\left| \alpha  \right|^2 ,\left| {h_j (I_\varepsilon  ) - h_j (I_{ex} )} \right| \leqslant C_3 r\left| \alpha  \right|^2 \varepsilon  \hfill \\
\end{gathered}
\]
for all $j\in\{1,2\}$, where $C_3$ is a positive constant depending only on $I_{ex}$. Therefore,
\[
\left| {g_j (I_{ex} )} \right| \leqslant \frac{{\left| {\alpha_j h_0 (I_{ex} )} \right|}}
{{\left| \alpha  \right|^2 }}.\left| {D_2 (I_{ex} )} \right| + \frac{{\left| {h_j (I_{ex} )} \right|}}
{{\left| \alpha  \right|^2 }}.\left| {D_1 (I_{ex} )} \right| \leqslant 2C_3 \left\| {\varphi _{ex} } \right\|_{L^1 (0,T)} r
\]
and
\[
\begin{gathered}
  \left| {g_j (I_\varepsilon  ) - g_j (I_{ex} )} \right| \leqslant \frac{{\left| {\alpha_j (h_0 (I_\varepsilon  ) - h_0 (I_{ex} ))} \right|}}
{{\left| \alpha  \right|^2 }}.\left| {D_2 (I_\varepsilon  )} \right| + \frac{{\left| {\alpha_j h_0 (I_{ex} )} \right|}}
{{\left| \alpha  \right|^2 }}.\left| {D_2 (I_\varepsilon  ) - D_2 (I_{ex} )} \right| \hfill \\
\h\h\h\h\h   + \frac{{\left| {h_j (I_\varepsilon  ) - h_0 (I_{ex} )} \right|}}
{{\left| \alpha  \right|^2 }}.\left| {D_1 (I_\varepsilon  )} \right| + \frac{{\left| {h_j (I_{ex} )} \right|}}
{{\left| \alpha  \right|^2 }}.\left| {D_1 (I_\varepsilon  ) - D_1 (I_{ex} )} \right| \hfill \\
   \leqslant C_3 r\varepsilon .\left( {\left\| {\varphi _{ex} } \right\|_{L^1 (0,T)}^2  + \varepsilon } \right) + C_3 r.\varepsilon  + C_3 r\varepsilon .\left( {\left\| {\varphi _{ex} } \right\|_{L^1 (0,T)}^2  + \varepsilon } \right) + C_3 r.\varepsilon  \hfill \\
   \leqslant 2C_3 \left( {\left\| {\varphi _{ex} } \right\|_{L^1 (0,T)}^2  + e^{ - e}  + 1} \right)r\varepsilon  \hfill \\
\end{gathered}
\]
\h Returning Lemma 4, for all $j\in\{1,2\}$, we get $G_j(I_\varepsilon)\in L^1(R^2)\cap L^2(R^2)$  because the support of $G_j(I_\varepsilon)$ is contained in $\overline {B(0,R_\varepsilon  )}$ and $G_j(I_\varepsilon)\in L^\infty  (R^2)$. Moreover,
\[
\begin{gathered}
  \left| {G_j (I_\varepsilon  ) - F(\widetilde f_{jex} )} \right| \leqslant \chi (B(0,R_\varepsilon  ))\left| {\frac{{g_j \left( {I_\varepsilon  } \right)D(I_\varepsilon  )}}
{{\delta _\varepsilon   + \left( {D(I_\varepsilon  )} \right)^2 }} - \frac{{g_j \left( {I_{ex} } \right)D(I_{ex} )}}
{{\delta _\varepsilon   + \left( {D(I_{ex} )} \right)^2 }}} \right| \hfill \\
  \h + \chi (B(0,R_\varepsilon  ))\left| {\frac{{g_j \left( {I_{ex} } \right)D(I_{ex} )}}
{{\delta _\varepsilon   + \left( {D(I_{ex} )} \right)^2 }} - {\frac{{g_j \left( {I_{ex} } \right)}}
{{D(I_{ex} )}}}} \right| + \chi (R^2 \backslash B(0,R_\varepsilon  )).\left| {F(\widetilde f_{jex} )} \right| \hfill \\
\end{gathered}
\]
We shall estimate each of the terms of the right-hand side. We have
\[
\begin{gathered}
  \left| {\frac{{g_j \left( {I_\varepsilon  } \right)D(I_\varepsilon  )}}
{{\delta _\varepsilon   + \left( {D(I_\varepsilon  )} \right)^2 }} - \frac{{g_j \left( {I_{ex} } \right)D(I_{ex} )}}
{{\delta _\varepsilon   + \left( {D(I_{ex} )} \right)^2 }}} \right| \leqslant \frac{{\delta _\varepsilon  \left| {g_j \left( {I_\varepsilon  } \right)D(I_\varepsilon  ) - g_j \left( {I_{ex} } \right)D(I_{ex} )} \right|}}
{{\left( {\delta _\varepsilon   + \left( {D(I_\varepsilon  )} \right)^2 } \right)\left( {\delta _\varepsilon   + \left( {D(I_{ex} )} \right)^2 } \right)}} \hfill \\
   + \frac{{\left| {D(I_\varepsilon  )} \right|.\left| {D(I_{ex} )} \right|.\left| {g_j \left( {I_\varepsilon  } \right)D(I_{ex} ) - g_j \left( {I_{ex} } \right)D(I_\varepsilon  )} \right|}}
{{\left( {\delta _\varepsilon   + \left( {D(I_\varepsilon  )} \right)^2 } \right)\left( {\delta _\varepsilon   + \left( {D(I_{ex} )} \right)^2 } \right)}} \hfill \\
   \leqslant \frac{{\left| {g_j \left( {I_\varepsilon  } \right)D(I_\varepsilon  ) - g_j \left( {I_{ex} } \right)D(I_{ex} )} \right|}}
{{\delta _\varepsilon  }} + \frac{{\left| {g_j \left( {I_\varepsilon  } \right)D(I_{ex} ) - g_j \left( {I_{ex} } \right)D(I_\varepsilon  )} \right|}}
{{\delta _\varepsilon  }} \hfill \\
\end{gathered}
\]
If $\varepsilon\in (0,e^{-e})$ then $R_{\varepsilon}>r_0$, so for all $\alpha\in B(0,R_\varepsilon)$ we get
\[
\begin{gathered}
  \left| {g_{j} (I_\varepsilon  )D (I_\varepsilon  ) - g_{j} (I_{ex} )D (I_{ex} )} \right| \hfill \\
   \leqslant \left| {g_{j} (I_\varepsilon  ) - g_{j} (I_{ex} )} \right|.\left| {D (I_\varepsilon  )} \right| + \left| {g_{j} (I_{ex} )} \right|.\left| {D (I_\varepsilon  ) - D (I_{ex} )} \right| \hfill \\
   \leqslant C_2R_\varepsilon\varepsilon .(C_2 + \varepsilon) + C_2R_\varepsilon\varepsilon  \leqslant (C_2+1)^2R_\varepsilon\varepsilon  \hfill \\
\end{gathered}
\]
and similarly,
\[
  \left| {g_{j} (I_\varepsilon  )D (I_{ex}  ) - g_{j} (I_{ex} )D(I_{\varepsilon} )} \right|\le (C_2+1)^2R_\varepsilon\varepsilon
\]
Consequently, for all $\varepsilon\in (0,e^{-e})$, we can estimate the first term
\[
\chi (B(0,R_\varepsilon  ))\left| {\frac{{g_j \left( {I_\varepsilon  } \right)D(I_\varepsilon  )}}
{{\delta _\varepsilon   + \left( {D(I_\varepsilon  )} \right)^2 }} - \frac{{g_j \left( {I_{ex} } \right)D(I_{ex} )}}
{{\delta _\varepsilon   + \left( {D(I_{ex} )} \right)^2 }}} \right| \leqslant \chi (B(0,R_\varepsilon  )).\frac{{2(C_2  + 1)^2 R_\varepsilon  \varepsilon }}
{{\delta _\varepsilon  }}
\]
Considering the second term, we have
\[
\left| {\frac{{g_j \left( {I_{ex} } \right)D(I_{ex} )}}
{{\delta _\varepsilon   + \left( {D(I_{ex} )} \right)^2 }} - \frac{{g_j \left( {I_{ex} } \right)}}
{{D(I_{ex} )}}} \right| = \frac{{\delta _\varepsilon  \left| {g_j \left( {I_{ex} } \right)} \right|}}
{{\left( {\delta _\varepsilon   + \left( {D(I_{ex} )} \right)^2 } \right).\left| {D(I_{ex} )} \right|}}
\]
We always have
\[
\frac{{\delta _\varepsilon  \left| {g_j \left( {I_{ex} } \right)} \right|}}
{{\left( {\delta _\varepsilon   + \left( {D(I_{ex} )} \right)^2 } \right).\left| {D(I_{ex} )} \right|}} \leqslant \left| {\frac{{g_j \left( {I_{ex} } \right)}}
{{D(I_{ex} )}}} \right| = 2\left| {\int\limits_\Omega  {f_{jex} (x)\cos (\alpha  \cdot x)dx} } \right| \leqslant 2\left\| {f_{jex} } \right\|_{L^2 (\Omega )}
\]
Furthermore, if $\alpha  \in B(0,R_\varepsilon  )\backslash B_\varepsilon $ then
\[
\frac{{\delta _\varepsilon  \left| {g_j \left( {I_{ex} } \right)} \right|}}
{{\left( {\delta _\varepsilon   + \left( {D(I_{ex} )} \right)^2 } \right).\left| {D(I_{ex} )} \right|}} \leqslant \frac{{\delta _\varepsilon  \left| {g_j \left( {I_{ex} } \right)} \right|}}
{{\left| {D(I_{ex} )} \right|^3 }} \leqslant \frac{{\delta _\varepsilon  C_2 R_\varepsilon  }}
{{\varepsilon ^{6q} }}
\]
Therefore, for all $\varepsilon\in (0,e^{-e})$, we can estimate the second term
\[
\chi (B(0,R_\varepsilon  ))\left| {\frac{{g_j \left( {I_{ex} } \right)D(I_{ex} )}}
{{\delta _\varepsilon   + \left( {D(I_{ex} )} \right)^2 }} - \frac{{g_j \left( {I_{ex} } \right)}}
{{D(I_{ex} )}}} \right| \leqslant 2\chi (B_\varepsilon  )\left\| {f_{jex} } \right\|_{L^2 (\Omega )}  + \chi (B(0,R_\varepsilon  ))\frac{{\delta _\varepsilon  C_2 R_\varepsilon  }}
{{\varepsilon ^{6q} }}
\]
Thus, for all $\varepsilon\in (0,e^{-e})$, we have
\[
\begin{gathered}
  \left| {G_j (I_\varepsilon  ) - F(\widetilde f_{jex} )} \right| \leqslant \chi (B(0,R_\varepsilon  ))\left( {\frac{{2(C_2  + 1)^2 R_\varepsilon  \varepsilon }}
{{\delta _\varepsilon  }} + \frac{{\delta _\varepsilon  C_2 R_\varepsilon  }}
{{\varepsilon ^{6q} }}} \right) \hfill \\
 \h\h  + 2\chi (B_\varepsilon  )\left\| {f_{jex} } \right\|_{L^2 (\Omega )}  + \chi (R^2 \backslash B(0,R_\varepsilon  ))\left| {F(\widetilde f_{jex} )} \right| \hfill \\
\end{gathered}
\]
Choosing $\delta _\varepsilon   = \varepsilon ^{\frac{{6q + 1}}{2}}$ and $C_0=2(C_2+1)^2+C_2$, we complete the proof.
\end{proof}
It is obvious that, for all $j\in\{1,2\}$, by Lebesgue's dominated
convergence theorem, $\chi (R^2 \backslash B(0,R_\varepsilon
))\left| {F(\widetilde f_{jex} )} \right|$ converges to 0 in
$L^2(R^2)$ when $\varepsilon\to 0$. However, to get an explicitly
estimate for it, some a-priori information about $f_{ex}$ must be
assume.
\begin{lemma} Let $a\in R$, $Q$ be an measurable subset of $R^n$ ($n\ge 1$), and $w\in L^1(Q)\cap L^2(Q)$. Then
\[
\int\limits_{R^n } {\left| {\int\limits_Q {w(x)\sin (a + \sum\limits_{k = 1}^n {\alpha _k x_k } )dx} } \right|^2 d\alpha }  = 2^{n - 1} \pi ^n \left\| w \right\|_{_{L^2 (Q)} }^2
\]
\end{lemma}
\begin{proof} We first prove in the case $a=0$. Put $\widetilde w:R^n\to R$
\[
\widetilde w (x)=\chi(Q)w(x)-\chi(-Q)w(-x)
\]
Then $\widetilde w \in L^1(R^n)\cap L^2(R^n)$ and
\[
F_n (\widetilde w)(\alpha ) = 2i\int\limits_Q {w(x)\sin (\sum\limits_{k = 1}^n {\alpha _k x_k } )dx}
\]
where $F_n$ is the Fourier transform in $R^n$. Using Paserval equality, we get
\[
\int\limits_{R^n } {\left| {\int\limits_Q {w(x)\sin (\sum\limits_{k = 1}^n {\alpha _k x_k } )dx} } \right|^2 d\alpha }  = \frac{1}
{4}\left\| {F_n (\widetilde w)} \right\|_{L^2 (R^n)}^2  = \frac{{(2\pi )^n }}
{4}\left\| {\widetilde w} \right\|_{L^2 (R^n)}^2  = 2^{n - 1} \pi ^n \left\| w \right\|_{L^2 (Q)}^2
\]
Similarly, we also have
\[
\int\limits_{R^n } {\left| {\int\limits_Q {w(x)\cos (\sum\limits_{k = 1}^n {\alpha _k x_k } )dx} } \right|^2 d\alpha }  = 2^{n - 1} \pi ^n \left\| w \right\|_{L^2 (Q)}^2
\]
Now, we notice that
\[
\begin{gathered}
  \left| {\int\limits_Q {w(x)\sin (a + \sum\limits_{k = 1}^n {\alpha _k x_k } )dx} } \right|^2  = (\cos (a))^2 \left| {\int\limits_Q {w(x)\sin (\sum\limits_{k = 1}^n {\alpha _k x_k } )dx} } \right|^2  \hfill \\
 \h\h  + (\sin (a))^2 \left| {\int\limits_Q {w(x)\cos (\sum\limits_{k = 1}^n {\alpha _k x_k } )dx} } \right|^2  + v(\alpha) \hfill \\
\end{gathered}
\]
where
\[
v(\alpha ) = \sin (2a)
.\int\limits_\Omega  {w(x)\sin (\sum\limits_{k = 1}^n {\alpha _k x_k } )dx} .\int\limits_\Omega  {w(x)\cos (\sum\limits_{k = 1}^n {\alpha _k x_k } )dx}
\]
Since $v(-\alpha)=-v(\alpha)$ for all $\alpha \in R^n$, we get $\int\limits_{R^n } {v(\alpha )d\alpha }  = 0$. Thus
\[
\begin{gathered}
  \int\limits_{R^n } {\left| {\int\limits_Q {w(x)\sin (a + \sum\limits_{k = 1}^n {\alpha _k x_k } )dx} } \right|^2 d\alpha }  \hfill \\
   = (\cos (a))^2 .2^{n - 1} \pi ^n \left\| w \right\|_{L^2 (Q)}^2  + (\sin (a))^2 .2^{n - 1} \pi ^n \left\| w \right\|_{L^2 (Q)}^2  = 2^{n - 1} \pi ^n \left\| w \right\|_{L^2 (Q)}^2  \hfill \\
\end{gathered}
\]
The proof is completed.
\end{proof}
Using Lemma 5, we have the following result.
\begin{lemma} Let $w\in H^1(\Omega)$ and $r>\pi/(2\sqrt{2})$. Then
$$
  \int\limits_{R^2 \backslash B(0,r)} {\left| {\int\limits_\Omega  {w(x )\cos (\alpha \cdot x )dx } } \right|^2d\alpha }  \leqslant \frac{{72\sqrt{2}\pi}}{ r}\left\| w \right\|_{H^1(\Omega)}^2 \hfill \\
$$
\end{lemma}
\begin{proof} Since
\[
\int\limits_{R^2 \backslash B(0,r)} {\left| {\int\limits_\Omega  {w(x)\cos (\alpha  \cdot x)dx} } \right|^2 d\alpha }  \leqslant \sum\limits_{j = 1}^2 {\int\limits_{\left| {\alpha _j } \right| \geqslant r/\sqrt 2 } {\left| {\int\limits_\Omega  {w(x)\cos (\alpha  \cdot x)dx} } \right|^2 d\alpha } }
\]
, the proof will be completed if we show that, for all $j\in\{1,2\}$,
\[
\int\limits_{\left| {\alpha _j } \right| \geqslant r/\sqrt 2 } {\left| {\int\limits_\Omega  {w(x)\cos (\alpha  \cdot x)dx} } \right|^2 d\alpha }  \leqslant \frac{{24\sqrt 2 \pi }}
{r}\left( {\left\| w \right\|_{L^2 (\Omega )}^2  + 2\left\| {\frac{{\partial w}}
{{\partial x_j }}} \right\|_{L^2 (\Omega )}^2 } \right)
\]
We will prove for the case $j=1$, and the other cases are similar. We have
\[
\int\limits_\Omega  {w(x)\cos (\alpha  \cdot x)dx}  = {\int\limits_0^1 {\left[ {w(x)\frac{{\sin (\alpha  \cdot x)}}
{{\alpha _1 }}} \right]} _{x_1  = 0}^{x_1  = 1} dx_2}  - \int\limits_\Omega  {\frac{{\partial w}}
{{\partial x_1}}.\frac{{\sin (\alpha  \cdot x)}}
{{\alpha _1 }}dx}
\]
so
\[
\begin{gathered}
  \left| {\int\limits_\Omega  {w(x)\cos (\alpha  \cdot x)dx} } \right|^2  \leqslant \frac{3}
{{\alpha _1^2 }}\left| { {\int\limits_0^1 {w(1,x_2)\sin (\alpha_1+\alpha _2 x_2)} dx_2 } } \right|^2  \hfill \\
   + \frac{3}
{{\alpha _1^2 }}\left| {{\int\limits_0^1 {w(0,x_2)\sin (\alpha _2 x_2)} dx_2} } \right|^2  + \frac{3}
{{\alpha _1^2 }}\left| {\int\limits_\Omega  {\frac{{\partial w}}
{{\partial x_1 }}.\sin (\alpha  \cdot x)dx} } \right|^2  \hfill \\
\end{gathered}
\]
Therefore,
\[
\begin{gathered}
  \int\limits_{\left| {\alpha _1 } \right| \geqslant r/\sqrt 2 } {\left| {\int\limits_\Omega  {w(x)\cos (\alpha  \cdot x)dx} } \right|^2 d\alpha }  \leqslant \frac{6}
{{r^2 }}\int\limits_{R^2 } {\left| {\int\limits_\Omega  {\frac{{\partial w}}
{{\partial x_1 }}(x).\sin (\alpha  \cdot x)dx} } \right|^2 d\alpha }  \hfill \\
   + \int\limits_{\left| {\alpha _1 } \right| \geqslant r/\sqrt 2 } {\frac{3}
{{\alpha _1^2 }}d\alpha _1 .} \int\limits_{ - \infty }^\infty  {\left| {\int\limits_0^1 {w(1,x_2 )\sin (\alpha _1  + \alpha _2 x_2 )dx_2 } } \right|^2 d\alpha _2 }  \hfill \\
   + \int\limits_{\left| {\alpha _1 } \right| \geqslant r/\sqrt 2 } {\frac{3}
{{\alpha _1^2 }}d\alpha _1 .} \int\limits_{ - \infty }^\infty  {\left| {\int\limits_0^1 {w(0,x_2 )\sin (\alpha _2 x_2 )dx_2 } } \right|^2 d\alpha _2 }  \hfill \\
   = \frac{{12\pi^2}}
{{r^2 }}\left\| {\frac{{\partial w}}
{{\partial x_1 }}} \right\|^2_{L^2 (\Omega )}  + \frac{{6\sqrt 2 \pi }}
{r}\left\| {w(1,.)} \right\|_{L^2 (0,1)}^2  + \frac{{6\sqrt 2 \pi  }}
{r}\left\| {w(0,.)} \right\|_{L^2 (0,1)}^2  \hfill \\
\end{gathered}
\]
Noting that
\[
w(1,x_2) = \int\limits_0^1 {\frac{\partial }
{{\partial x_1 }}\left( {x_1 w(x)} \right)dx_1 }  = \int\limits_0^1 {\left( {w(x) + x_1 \frac{{\partial w}}
{{\partial x_1 }}(x)} \right)dx_1 }
\]
, we get
\[
\left| {w(1,x_2)} \right|^2  \leqslant \int\limits_0^1 {\left( {2\left| {w(x)} \right|^2  + 2\left| {\frac{{\partial w}}
{{\partial x_1 }}(x)} \right|^2 } \right)dx_1 }
\]
Hence,
\[
{\int\limits_0^1 {\left| {w(1,x_2 )} \right|^2 dx_2 } }  \leqslant 2\left\| w \right\|_{L^2 (\Omega )}^2  + 2\left\| {\frac{{\partial w}}
{{\partial x_1 }}} \right\|_{L^2 (\Omega )}^2
\]
Similarly,
\[
\begin{gathered}
  {\int\limits_0^1 {\left| {w(0,x_2 )} \right|^2 dx_2 } }  =  {\int\limits_0^1 {\left| {\int\limits_0^1 {\frac{\partial }
{{\partial x_1 }}\left( {(1 - x_1 )w(x)} \right)dx_1 } } \right|^2 dx_2 } }  \hfill \\
   \leqslant  {\int\limits_0^1 {\int\limits_0^1 {\left( {2\left| {w(x)} \right|^2  + 2\left| {\frac{{\partial w}}
{{\partial x_1 }}(x)} \right|^2 } \right)dx_1 } dx_2 } }= 2\left\| w \right\|_{L^2 (\Omega )}^2  + 2\left\| {\frac{{\partial w}}
{{\partial x_1 }}} \right\|_{L^2 (\Omega )}^2  \hfill \\
\end{gathered}
\]
Thus, we have
\[
\begin{gathered}
  \int\limits_{\left| {\alpha _1 } \right| \geqslant r/\sqrt 2 } {\left| {\int\limits_\Omega  {w(x)\cos (\alpha  \cdot x)dx} } \right|^2 d\alpha }  \leqslant \frac{{12\pi ^2 }}
{{r^2 }}\left\| {\frac{{\partial w}}
{{\partial x_1 }}} \right\|_{L^2 (\Omega )}  +  \hfill \\
   + \frac{{24\sqrt 2 \pi }}
{r}\left( {\left\| {w(1,.)} \right\|_{L^2 (\Omega )}^2  + \left\| {\frac{{\partial w}}
{{\partial x_1 }}} \right\|_{L^2 (\Omega )}^2 } \right) \leqslant \frac{{24\sqrt 2 \pi }}
{r}\left( {\left\| {w(1,.)} \right\|_{L^2 (\Omega )}^2  + 2\left\| {\frac{{\partial w}}
{{\partial x_1 }}} \right\|_{L^2 (\Omega )}^2 } \right) \hfill \\
\end{gathered}
\]
The proof is completed.
\end{proof}
\begin{remark} By the same way, we can show that, if $w\in H^1(\Omega)$ and $r>\pi/(2\sqrt{2})$ then
$$
\int\limits_{R^2 \backslash B(0,r)} {\left| {\int\limits_Q {w(x_1 ,x_2 )\cos (\alpha_1 x_1 )\cos (\alpha_2 x_2 )dx } } \right|^2d\alpha}  \leqslant \frac{{16\sqrt{2}\pi}}{ r}\left\| w \right\|_{H^1(Q)}^2
$$
This result improves immediately the results of \cite{TQD}.
\end{remark}
\text{}\\{\bf Proof of theorem 2}
\begin{proof} Recall that $q, \delta _\varepsilon, R_\varepsilon$ are defined by $(\ref{ths})$, and $G_j(I_\varepsilon)$, $B_{\varepsilon}$ are defined by Lemma 4. For all $j\in\{1,2\}$, we define $ f_{j\varepsilon}:R^2\to R$
\[
 f_{j\varepsilon } (\xi) =\frac{1}{4\pi^2} \int\limits_{R^2 } {G_j (I_\varepsilon  )(\alpha)e^{i(\xi \cdot \alpha)} d\alpha}
\]
Applying Lemma 4, we have $G_j(I_\varepsilon)\in L^1(R^2)\cap
L^2(R^2)$ , so $f_{j\varepsilon }\in C(R^2)\cap L^2(R^2)$ and $F(
f_{j\varepsilon})=G_j(I_\varepsilon)$. Applying Lemma 4 again, for
all $\varepsilon \in (0,e^{-e})$, we get \bq
\begin{gathered}
  \left| {F(f_{j\varepsilon } ) - F(\widetilde f_{jex} )} \right| \leqslant \chi (B(0,R_\varepsilon  ))C_0 R_\varepsilon  \varepsilon ^{\frac{{1 - 6q}}
{2}}  \hfill \\
  \h + 2\chi (B_\varepsilon  )\left\| {f_{jex} } \right\|_{L^2 (\Omega )}  + \chi (R^2 \backslash B(0,R_\varepsilon  ))\left| {F(\widetilde f_{jex} )} \right| \hfill \\
\end{gathered}
\label{est}
\eq
where $C_0$ is a positive constant depending only $I_{ex}$. It implies that
\[
\begin{gathered}
  \left| {F(f_{j\varepsilon}  ) - F(\widetilde f_{jex} )} \right|^2  \leqslant 2\chi (B(0,R_\varepsilon  ))C_0^2 R_\varepsilon ^2 \varepsilon ^{1 - 6q}  \hfill \\
 \h  + 4\chi (B_\varepsilon  )\left\| {f_{jex} } \right\|_{L^2 (\Omega )}^2  + 2\chi (R^2 \backslash B(0,R_\varepsilon  ))\left| {F(\widetilde f_{jex} )} \right|^2  \hfill \\
\end{gathered}
\]
Hence,
\[
\left\| {F(f_{j\varepsilon } ) - F(\widetilde f_{jex} )} \right\|_{L^2 (R^2)}^2  \leqslant 2C_0^2 \pi R_\varepsilon ^4 \varepsilon ^{1 - 6q}  + 4m(B_\varepsilon  )\left\| {f_{jex} } \right\|_{L^2 (\Omega )}^2  + 2\int\limits_{R^2 \backslash B(0,R_\varepsilon  )} {\left| {F(\widetilde f_{jex} )} \right|^2 d\alpha }
\]
It is obvious that $2C_0^2 \pi R_\varepsilon ^4 \varepsilon ^{1 - 6q} \leqslant R_\varepsilon ^{-1}$ for $\varepsilon>0$ small enough. Moreover, since
\[
B_\varepsilon   \subset \left( {\left\{ {\alpha  \in B(0,R_\varepsilon  ),\left| {D_1 (I_{ex} )(\alpha )} \right| \leqslant \varepsilon ^q } \right\} \cup \left\{ {\alpha  \in B(0,R_\varepsilon  ),\left| {D_2 (I_{ex} )(\alpha )} \right| \leqslant \varepsilon ^q } \right\}} \right)
\]
, we apply Lemma 2 (with $\varphi_0(t)=\varphi_{ex}(T-t)$) to get that $m( B_\varepsilon)\le 2R_\varepsilon^{-1}$ for $\varepsilon>0$ small enough. Thus, for $\varepsilon>0$ small enough, we get
\[
\left\| {F(f_{j\varepsilon } ) - F(\widetilde f_{jex} )} \right\|_{L^2 (R^2 )}^2  \leqslant \frac{1}
{{R_\varepsilon  }} + \frac{8}
{{R_\varepsilon  }}\left\| {f_{jex} } \right\|_{L^2 (\Omega )}^2  + 2\int\limits_{R^2 \backslash B(0,R_\varepsilon  )} {\left| {F(\widetilde f_{jex} )} \right|^2 d\alpha d\beta }
\]
\h By Parseval equality, we have \bq
\begin{gathered}
  \left\| {f_{j\varepsilon }  - f_{jex} } \right\|_{L^2 (\Omega )}^2  \leqslant \left\| {f_{j\varepsilon }  - \widetilde f_{jex} } \right\|_{L^2 (R^2 )}^2  = \frac{1}
{{4\pi ^2 }}\left\| {F(f_{j\varepsilon } ) - F(\widetilde f_{jex} )} \right\|_{L^2 (R^2 )}^2  \hfill \\
   \leqslant \frac{1}
{{4\pi ^2 }}\left( {\frac{1}
{{R_\varepsilon  }} + \frac{8}
{{R_\varepsilon  }}\left\| {f_{jex} } \right\|_{L^2 (\Omega )}^2  + 2\int\limits_{R^2 \backslash B(0,R_\varepsilon  )} {\left| {F(\widetilde f_{jex} )} \right|^2 } d\alpha } \right) \hfill \\
\end{gathered}
\label{dgss}
\eq
for $\varepsilon>0$ small enough. Since $F(\widetilde f_{jex})\in L^2(R^2)$, we obtain that
\[
\mathop {\lim }\limits_{\varepsilon  \to 0} \left\| {f_{j\varepsilon }  - f_{jex} } \right\|_{L^2 (\Omega )}  = 0
\]
\h If $f_{jex}\in H^1(\Omega)$ then using $(\ref{dgss})$ and Lemma 6, we get
\[
\begin{gathered}
  \left\| {f_{j\varepsilon }  - f_{jex} } \right\|_{L^2 (\Omega )}^2  \leqslant \frac{1}
{{4\pi ^2 }}\left( {\frac{1}
{{R_\varepsilon  }} + \frac{8}
{{R_\varepsilon  }}\left\| {f_{jex} } \right\|_{L^2 (\Omega )}^2  + 2.4.\frac{{72\sqrt2 \pi  }}
{{R_\varepsilon  }}\left\| {f_{jex} } \right\|_{H^1 (\Omega )}^2 } \right) \hfill \\
 \leqslant \left( {66\left\| {f_{jex} } \right\|_{H^1 (\Omega )}^2  + \frac{1}
{{4\pi ^2 }}} \right).\frac{1}
{{R_\varepsilon  }} = 63eT\left( {66\left\| {f_{jex} } \right\|_{H^1 (\Omega )}^2  + \frac{1}
{{4\pi ^2 }}} \right).\frac{{\ln (\ln (\varepsilon ^{ - 1} ))}}
{{\ln (\varepsilon ^{ - 1} )}}
 \hfill \\
\end{gathered}
\]
for $\varepsilon>0$ small enough. This complete the proof.
\end{proof}
\text{}\\{\bf Proof of Theorem 3}
\begin{proof} We shall use the notations of the proof of Theorem 2. Notice that the assumtion
\[
\int\limits_{R^2 } {\left| {\int\limits_\Omega  {f_{jex} (x).\cos (\alpha  \cdot x)dx} } \right|d\alpha  < \infty },
\]
is equivalent to $F(\widetilde f_{jex})\in L^1(R^2)$. Since $\widetilde f_{jex}, F(\widetilde f_{jex})\in L^1(R^2)\cap L^2(R^2)$, we get
\[
\widetilde f_{jex} (\xi ) = \frac{1}{4\pi^2}\int\limits_{R^2 } {F(\widetilde f_{jex} )(\alpha )e^{i(\alpha \cdot \xi)} d\alpha}
\]
Therefore,
\bq
4\pi ^2 \left\| {f_{j\varepsilon }  - f_{jex} } \right\|_{L^\infty  (\Omega )}  \leqslant 4\pi ^2 \left\| {f_{j\varepsilon }  - \widetilde f_{jex} } \right\|_{L^\infty  (R^2 )}  \leqslant \left\| {F(f_{j\varepsilon } ) - F(\widetilde f_{jex} )} \right\|_{L^{1}(R^2 )}
\label{beq}
\eq
From $(\ref{est})$, we have
\[
\left\| {F(f_{j\varepsilon } ) - F(\widetilde f_{jex} )} \right\|_{L^1 (R^2 )}  \leqslant C_0 \pi R_\varepsilon ^3 \varepsilon ^{\frac{{1 - 3q}}
{2}}  + 2m(B_\varepsilon  )\left\| {f_{jex} } \right\|_{L^2 (\Omega )}  + \int\limits_{R^2 \backslash B(0,R_\varepsilon  )} {\left| {F(\widetilde f_{jex} )} \right|d\alpha }
\]
For $\varepsilon>0$ small enough, we have $C_0 \pi R_\varepsilon ^3\varepsilon ^{\frac{{1 - 3q}}
{2}} \le R _{\varepsilon}^{-1}$ and $m(B_\varepsilon  )\le 2R _{\varepsilon}^{-1}$. Thus, from $(\ref{beq})$, for $\varepsilon>0$ small enough, for all $j\in\{1,2\}$, we get
\[
4\pi ^2 \left\| {f_{j\varepsilon }  - f_{jex} } \right\|_{L^\infty  (\Omega )}  \leqslant \frac{1}
{{R_\varepsilon  }} + \frac{4}
{{R_\varepsilon  }}\left\| {f_{jex} )} \right\|_{L^2 (\Omega )}  + \int\limits_{R^2\backslash B(0,R_\varepsilon  )} {\left| {F(\widetilde f_{jex} )} \right|d\alpha }
\]
Since $F(\widetilde f_{jex} )\in L^1(R^2)$, we obtain that $\mathop {\lim }\limits_{\varepsilon  \to 0} \left\| {f_{j\varepsilon }  - f_{jex} } \right\|_{L^\infty  (\Omega )}  = 0$ for all $j\in\{1,2\}$.
\end{proof}

\text{}\\
{\bf Remark 2.} {\it We can replace $R_{\varepsilon}$ defined by
$(7)$ by
\[
\widetilde R_\varepsilon   = 10\left( {\ln (\varepsilon ^{ - 1} )}
\right)^{9/10}
\]
to construct a better regularized solution in the case that
$\varepsilon$ is not too small.}
\\\\
{\bf 4. A numerical experience}\\\\
Assume that $T=1$, $\mu=1/12$, $\lambda=-1/8$.
\\We consider the exact data $I_{ex}=(\varphi,X,u_0,u_0^*,u_T)$ given by
\[
\begin{gathered}
  \varphi  = \frac{{\pi ^2 }}
{3}\sin (\pi t), \hfill \\
  X_1  = \frac{\pi }
{6}\sin (\pi t).\left[ {\sin (2\pi x_2 )n_1  + \sin (4\pi x_1 )n_2 } \right], \hfill \\
  X_2  = \frac{\pi }
{6}\sin (\pi t).\left[ {\sin (2\pi x_1 )n_2  + \sin \left( {4\pi x_2 } \right)n_1 } \right], \hfill \\
  u_0  = u_T  = (0,0), \hfill \\
  u_0^*  = \left( {\pi \sin (4\pi x_1 )\sin (2\pi x_2 ),\pi \sin (2\pi x_1 )\sin (4\pi x_2 )} \right) .\hfill \\
\end{gathered}
\]
Then the corresponding exact solution of the system $(1)-(4)$ is
\[
\begin{gathered}
  u_{ex}  = \left( {\sin (\pi t)\sin (4\pi x_1 )\sin (2\pi x_2 ),\sin (\pi t)\sin (4\pi x_1 )\sin (2\pi x_2 )} \right), \hfill \\
  f_{ex}  = \left( {\cos (2\pi x_1 )\cos (4\pi x_2 ),\cos (4\pi x_1 )\cos (2\pi x_2 )} \right). \hfill \\
\end{gathered}
\]
For each $n=1,2,3,...$, we consider the inexact data
$I_{n}=(\varphi_n,X^n,u_0^n,u_0^{*n},u_T^n)$ given by
\[
\begin{gathered}
  \varphi _n  = \varphi , \hfill \\
  X_1^n  = X_1  + \frac{\pi }
{{12\sqrt n }}\sin (\pi t).\left[ {\sin (2n\pi x_2 )n_1  + 2\sin (2n\pi x_1 )n_2 } \right], \hfill \\
  X_2^n  = X_2  + \frac{\pi }
{{12\sqrt n }}\sin (\pi t).\left[ {\sin (2n\pi x_1 )n_2  + 2\sin (2n\pi x_2 )n_1 } \right], \hfill \\
  u_0^n  = u_T^n  = (0,0), \hfill \\
  u_0^{*n}  = u_0^*  + \frac{\pi }
{{n\sqrt n }}\sin (2n\pi x_1 )\sin (2n\pi x_2 )\left( {1,1} \right). \hfill \\
\end{gathered}
\]
Then the corresponding disturbed solution of the system $(1)-(4)$ is
\[
\begin{gathered}
 u^n  = u_{ex}  + \frac{1}
{{n\sqrt n }}\sin (\pi t)\sin (2n\pi x_1 )\sin (2n\pi x_2 )\left( {1,1} \right) ,\hfill \\
 f_{di}^n = f_{ex}  + \left[ {\left( {\frac{3}
{2}\sqrt n  - \frac{3} {{n\sqrt n }}} \right)\sin (2n\pi x_1 )\sin
(2n\pi x_2 ) + \frac{{\sqrt n }}
{2}\cos (2n\pi x_1 )\cos (2n\pi x_2 )} \right]\left( {1,1} \right). \hfill \\
\end{gathered}
\]
We get
\[
\begin{gathered}
  \varphi _n  = \varphi , \hfill \\
  \left\| {X_j^n (t,.) - X_j^{ex} (t,.)} \right\|_{L^1 (0,T,\partial \Omega )}  = \frac{2}
{{\pi \sqrt n }},\hfill\\
u_0^n  = u_0 ,u_T^n  = u_T ,\hfill\\
\left\| {u_{0j}^{*n}  - u_{0j}^* } \right\|_{L^1 (\Omega )}  =
\frac{4}
{{\pi n\sqrt n }},\forall j\in\{1,2\}, \hfill \\
\end{gathered}
\]
and
\[
\left\| { f_{jdi}^n  - f_{jex} } \right\|_{L^2 (\Omega )}^2  =
\frac{5} {8}n - \frac{9} {{4n}} + \frac{9} {{4n^3 }}.
\]
Hence, when $n$ is large, a small error of data will cause a large
error of solution. It show that the problem is ill$-$posed, and a
regularization is necessary.
\\\h We shall construct the regularized solution as in Theorem 1 corresponding $\varepsilon  = n^{ - 1/2}$. From the straightforward calculation, we obtain that
\[
\begin{gathered}
  D(I_n )(\alpha ) = \frac{{32\pi ^6 \sin \left( {\frac{{\left| \alpha  \right|}}
{{2\sqrt 6 }}} \right)\sin \left( {\frac{{\left| \alpha  \right|}}
{{2\sqrt 3 }}} \right)}}
{{\left( {\left| \alpha  \right|^2  - 24\pi ^2 } \right).\left( {\left| \alpha  \right|^2  - 12\pi ^2 } \right)}}, \hfill \\
  g_1 (I_n )(\alpha) = D(I_n )(\alpha ) \times (\sin (\alpha _1 )\sin (\alpha _2 ) - (1 - \cos (\alpha _1 ))(1 - \cos (\alpha _2 ))) \times  \hfill \\
 \h\h\h~\times \left( {\frac{{2\alpha _1 \alpha _2 }}
{{(\alpha _1^2  - 4\pi ^2 )(\alpha _2^2  - 16\pi ^2 )}} +
\frac{{\sqrt n (\alpha _1 \alpha _2  + 12\pi ^2 (2-n^2))}}
{{(\alpha _1^2  - 4n^2 \pi ^2 )(\alpha _2^2  - 4n^2 \pi ^2 )}}} \right). \hfill \\
\end{gathered}
\]
Thus, the regularized solution defined by
\[
f_{1re}^n (x) = \frac{1}{4\pi^2}\int\limits_{B(0,\widetilde R_n )}
{\frac{{g_1 (I_n )(\alpha ).D(I_n )(\alpha )}} {{\delta _n  + \left(
{D(I_n )(\alpha )} \right)^2 }}} .\cos (\alpha  \cdot x)d\alpha ,
\]
where
 \[
\delta _n  = n^{ - 13/28} ,\widetilde R_n  = 10\left( {\ln (\sqrt n
)} \right)^{9/10}.
\]
\h For example, if $\varepsilon=10^{-2}$ then
$$
\begin{gathered}
  n=10^4, \delta _n= 0.01389495494, \widetilde R_n= 39.52948133, \hfill \\
\end{gathered}
$$
and we have some figures about the exact solution $f_{1ex}$, the
disturbed solution $f_{1di}^n$ and the regularized
solution$f_{1re}^{n}$. \ \par
\ \par \vspace{1.5cm}
 \centerline{
\includegraphics[width=4.25in]{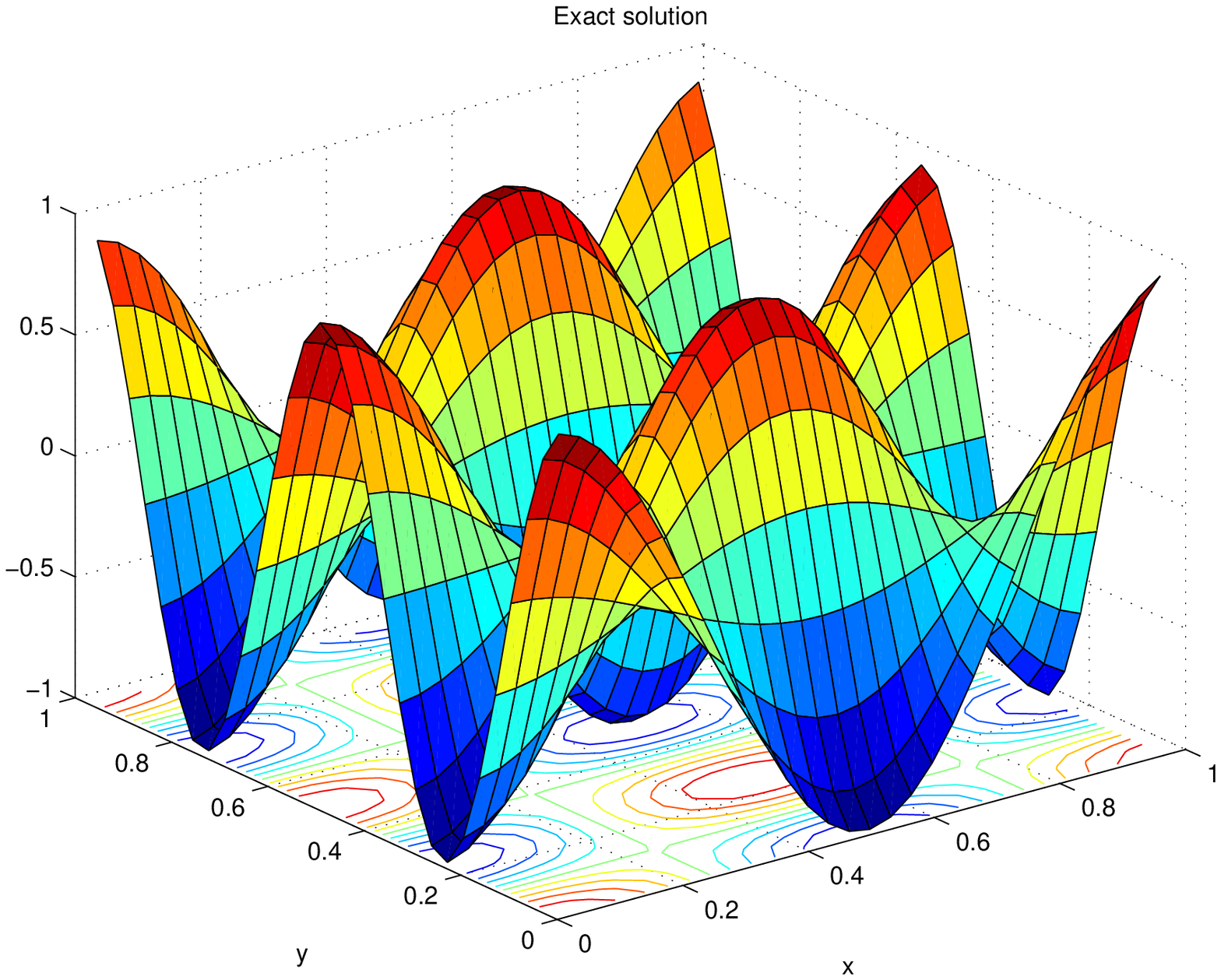}}
\center{Figure 1. The exact solution.}\endcenter
\ \par \centerline{
\includegraphics[width=4.25in]{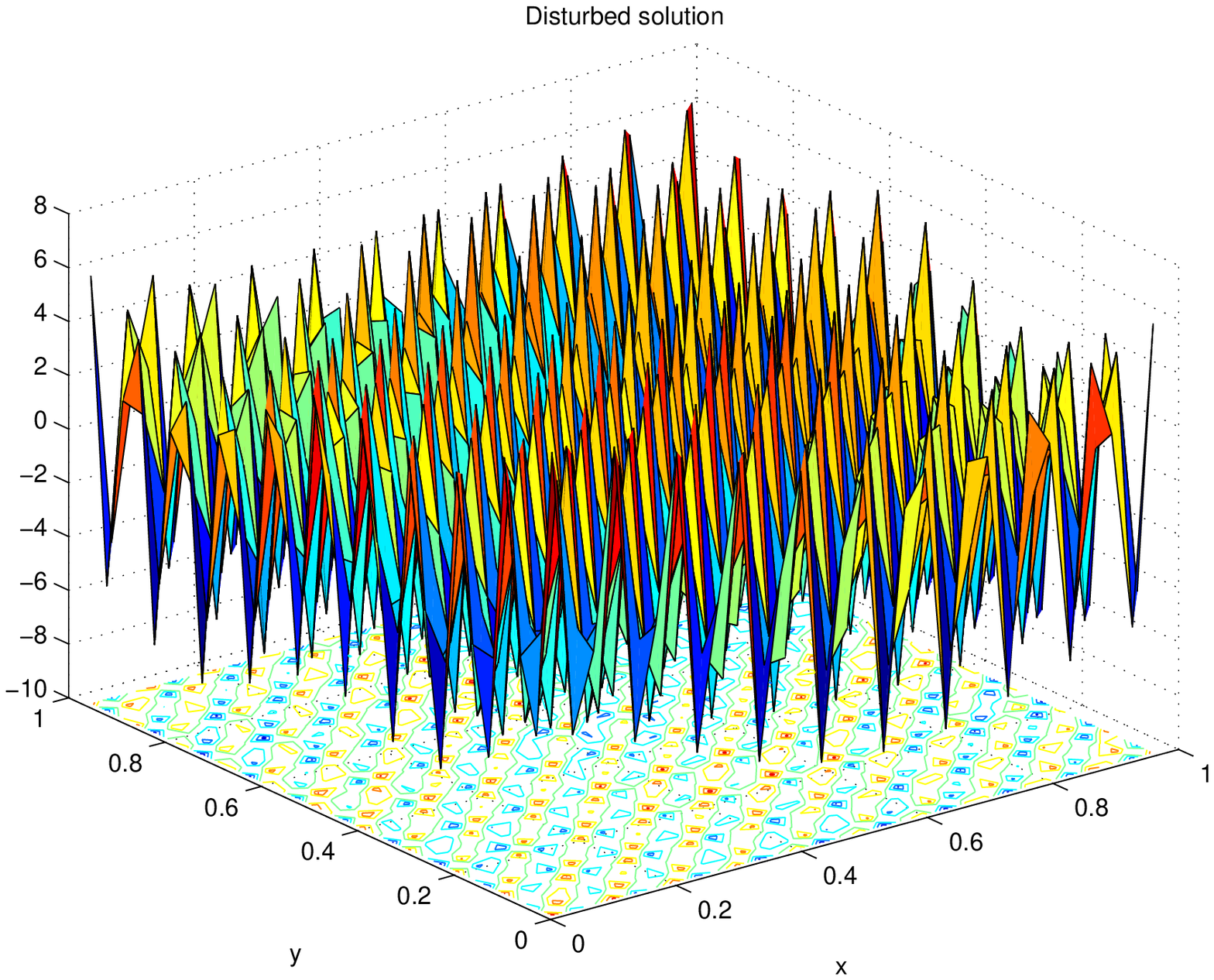}}
\center{Figure 2. The disturbed solution.}\endcenter
\ \par \centerline{
\includegraphics[width=4.5in]{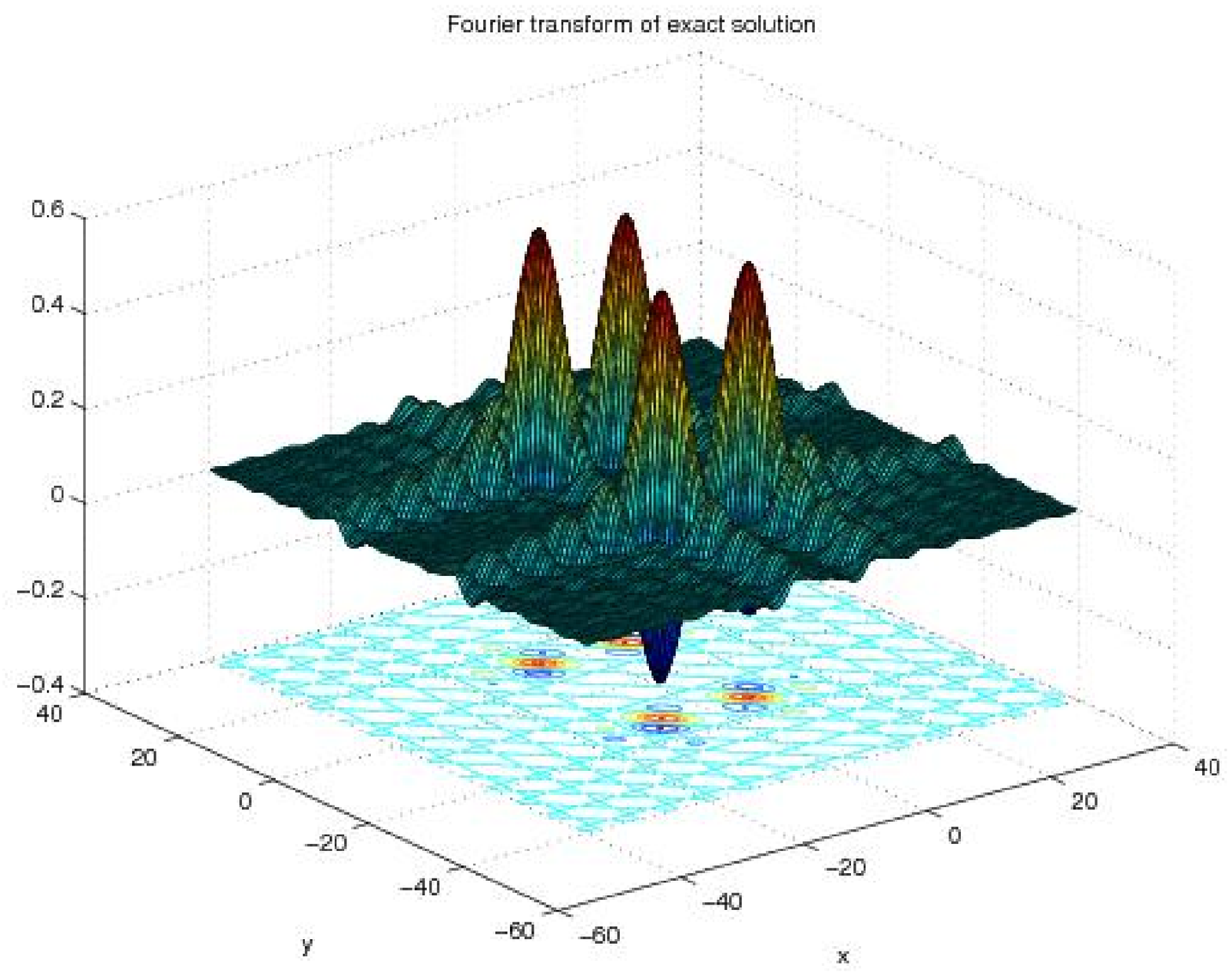}}
\center{Figure 3. The Fourier transform of the exact
solution.}\endcenter
\ \par \centerline{
\includegraphics[width=5in]{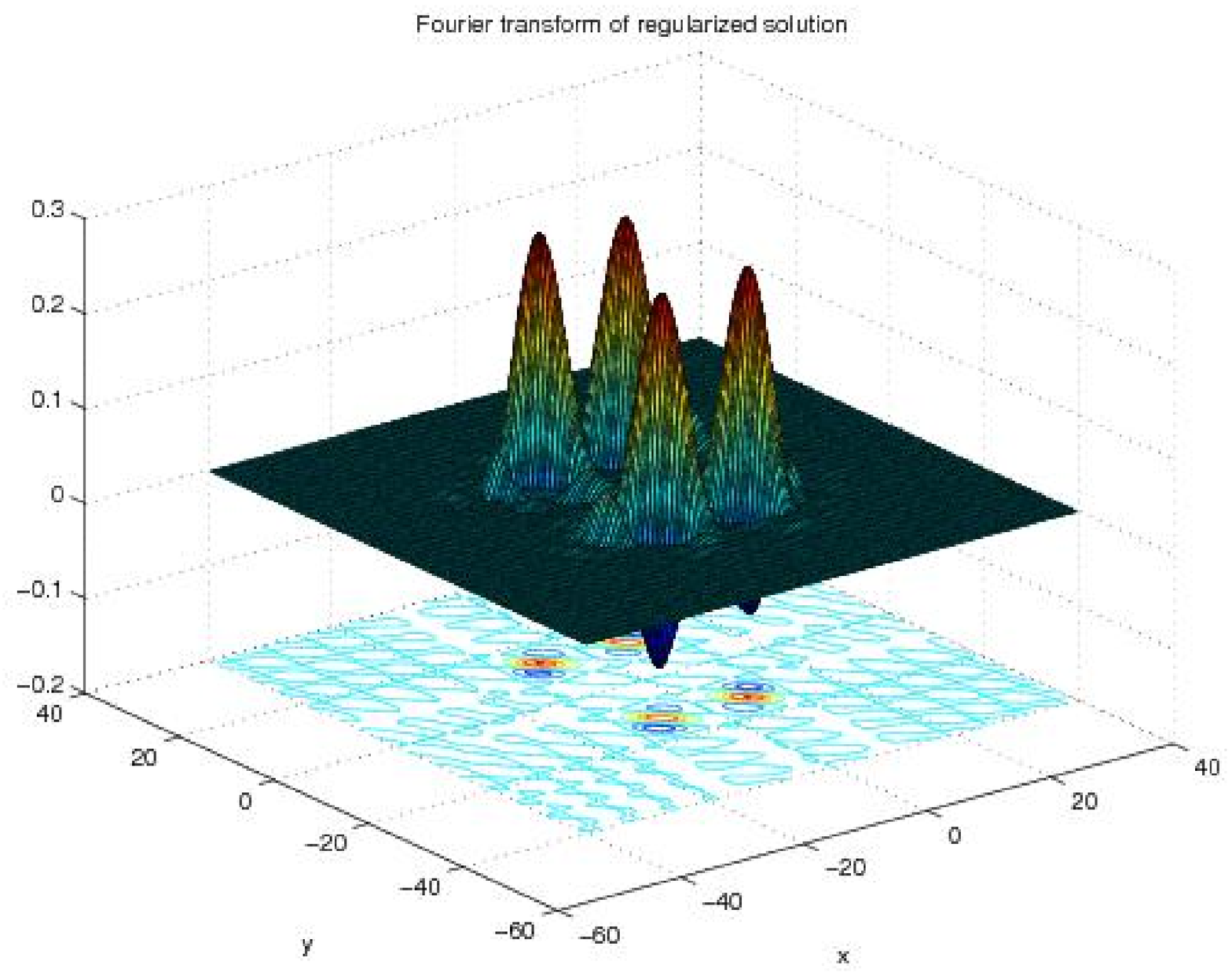}}
\center{Figure 4. The Fourier transform of the regularized
solution.}\endcenter

\end{document}